\documentclass{article}

\usepackage{amsmath, amssymb, amsthm, amsfonts, bm}
\usepackage{enumerate}

\numberwithin{equation}{section}
\newtheorem{cthm}{Theorem}
\newenvironment{ctheorem}[1]
  {\renewcommand\thecthm{#1}\cthm}
  {\endcthm}

\newtheorem{lemma}{Lemma}[section]
\theoremstyle{definition}

\theoremstyle{remark}

\newcommand{\Rnum}[1]{\uppercase\expandafter{\romannumeral #1\relax}}

\newcommand{\mc}[1]{\mathcal{#1}}

\title{A second alternative approach for the study of the Muckenhoupt class $A_1(\mathbb{ R})$}
\author{Eleftherios N. Nikolidakis}
\date{}

\begin{document}
\maketitle
\footnotetext{ {\em E-mail address}: lefteris@math.uoc.gr}

\begin{abstract}
We find the exact best possible range of those $p >1$ for which any $\varphi \in A_1(\mathbb{R})$, with $A_1$-constant equal to c, must also belong to $L^p$. In this way we provide alternative proofs of the results in \cite{2} and \cite{10}.
\end{abstract}

\section{Introduction} \label{sec:1}
The study of Muckenhoupt weights has been proved to be important in analysis. One of the most important facts about these is their self improving property. A way to express this is through the so-called reverse H\"{o}lder inequalities (see \cite{3}, \cite{4} and \cite{6}).

For an interval $\mc J$ on $\mathbb{R}$, we define the class $A_1(\mc J)$ to be the set of all those $\varphi: \mc J \rightarrow \mathbb{R}^+$ for which there exists a constant $c\geq1$, such that the following inequality is satisfied
\begin{equation} \label{eq:1p1}
\frac{1}{|\mc I|}\int_{\mc I} \varphi(x)dx \leq c\cdot essin f_{\mc I}(\varphi),
\end{equation}
for every subinterval $\mc I$ of $\mc J$, where $|^.|$ is the Lebesque measure on $\mathbb{R}$. The least constant $c$ for which \eqref{eq:1p1} holds, is called the $A_1$-constant of $\varphi$ and is denoted by $[\varphi]_1$. We will say then that $\varphi$ belongs to the class $A_1(\mc J)$ with constant c, and we will write $\varphi \in A_1(J,c)$.

The study of weights in the class $A_1(\mc J,c)$ has been seen for the first time in \cite{2}. In that paper the study of such weights has been given through the notion of the non-increasing rearrangement of $\varphi$, denoted by $\varphi^*$, which is a non-negative and non-increasing function defined on $(0,|\mc J|]$. It is characterized by the following two additional properties. It is equimeasurable to $\varphi$ (in the sense that $|\{\varphi > \lambda\}|=|\{\varphi^* >\lambda\}|$, for every $\lambda>0$) and is also left continuous. All these properties define uniquely $\varphi^*$ as can be seen in\cite{1}, \cite{5} or \cite{8}. Nevertheless an equivalent definition of $\varphi^*$ can be given by the following formula:
\begin{equation*}
\varphi^*(t)=\sup_{\substack{E\subseteq \mc J \\ |E|=t}}[\inf_{x\in E}\varphi(x)], \textrm{for t} {\in (0,|\mc J|]},
\end{equation*}
as can be seen in \cite{8}.

In \cite{2} now the following result has been proved:
\begin{ctheorem}{1} \label{thm:1}
Let $\varphi \in A_1(\mc J,c)$. The $\varphi^*$ satisfies:
\begin{equation} \label{eq:1p2}
\frac{1}{t} \int_{0}^{t} \varphi^*(y)dy \leq c \varphi^*(t), \text{for t} \in (0,|\mc J|].
\end{equation}
That is $\varphi^*$ belongs to the class $A_1(\mc J)$, with $A_1$-constant not more than c.
\end{ctheorem}
The above theorem describes the $A_1$-properties of $\varphi^*$, in terms of those of $\varphi$. It was used effectively by the authors in \cite{2} in order to  prove the following:
\begin{ctheorem}{2} \label{thm:2}
Let $\varphi \in A_1(\mc J,c)$. Then $\varphi \in L^p$ for every $p\in [1 ,\frac{c}{(c-1)})$. Moreover, the following inequality must hold for every subinterval $\mc I$ of $\mc J$, and every $p$ in the above range,
\begin{equation} \label{eq:1p3}
\frac{1}{|\mc I|}\int_{\mc I} \varphi^p(x)dx \leq \frac{1}{c^{p-1}(c+p-pc)}\Big(\frac{1}{|\mc I|}\int_{\mc I} \varphi(x)dx \Big)^p.
\end{equation}
Additionally, the above inequality is sharp, that is the constant appearing in the right side cannot be decreased.
\end{ctheorem}
Our aim in this paper is to give an alternative proof of Theorem 2, by using Theorem 1 and certain techniques involving the well known Hardy operator on the line. Additionally, we need to mention that in \cite{7} and \cite{9} related problems for estimates for the range of $p$ in higher dimensions have been treated.

The paper is organized as follows: In Section 2 we give a brief discussion of the proof of the Theorem 1, as is presented in \cite{2} and in Section 3 we provide the proof of Theorem 2.

\section{ $\varphi^*$ as an $A_1$ weight on $\mathbb{R}$.} \label{sec:2}
%$\textbf{Section 2}$: $\varphi^*$ as an $A_1$ weight on $\mc R$.
%\newline
Before we present the proof of Theorem 1 we give the following covering Lemma as can be seen in \cite{2}.
\begin{lemma} \label{lem:2p1}
 Let E be a measurable bounded subset of $\mathbb{R}$ and $\epsilon>0$. More precisely suppose that $E\subseteq I$, for a certain bounded interval $I$ of $\mathbb{R}$.  Then there exists a sequence $(I_\nu)_{\nu=1}^{\infty}$ of subintervals of $I$ with disjoint interiors and a subset $E_1$ of $E$ with the properties that $|E_1|=|E|$ and
\begin{enumerate}[i)]{\em
\item $E_1\subseteq \bigcup\limits_{r=1}^{\infty}I_{\nu} $
\item $(1-\epsilon) |I_{\nu}| \leq |I_{\nu}\cap E|<|I_{\nu}|$, for every $\nu$.
\vspace{-5pt}}
\end{enumerate}
\end{lemma}

We now proceed to the
\newline
$\underline{\textbf{Proof of Theorem 1}}:$ Suppose without loss of generality that $\mc J=(0,1)$ and that $\varphi$ satisfies \eqref{eq:1p1}, for every subinterval $\mc I$ on the above interval. Fix $t\in (0,1]$ and $\epsilon>0$. Let $E_t$ be a subset of $(0,1)$ such that $|E_t|=t$ and $\varphi(x)\leq \varphi^*(t)$, for any $x\notin E_t$. Using Lemma 2.1 we produce a subset $E_{t,1}$ of $E_t$, such that $|E_{t,1}|=t$ and $E_{t,1}\subseteq \bigcup\limits_{\nu=1}^{\infty}I_{\nu}$, where for every $\nu=1,2,...$ the following holds:
\begin{equation} \label{eq:2p1}
(1-\epsilon)|I_\nu|\leq |I_\nu \cap E_t|<|I_\nu|,
\end{equation}
for a suitable family $(I_\nu)_{\nu=1}^{\infty}$ of subintervals of $(0,1)$.
By the strict inequality in \eqref{eq:2p1}, we conclude that $I_\nu$ contains a set of positive measure in the complement of $E_t$, therefore we must have that
\begin{equation*}
essinf_{x\in I_\nu} \varphi(x)\leq \varphi^*(t),
\end{equation*}
so by using \eqref{eq:1p1} and \eqref{eq:2p1} we have a consequence that
\begin{equation*}
\begin{aligned}
& \int_{0}^{t} \varphi^*(y)dy=\int_{E_t} \varphi(x)dx=\int_{E_t,1} \varphi(x)dx\leq \sum_{\nu=1}^{\infty} \int_{I_\nu}\varphi(x)dx\leq c \sum_{\nu=1}^{\infty}|I_\nu|\cdot \varphi^*(t)\leq \\
&  \leq \frac{c}{1-\epsilon}\Big(\sum_{\nu=1}^{\infty}|I_\nu \cap E_t|\Big)\cdot \varphi^*(t)=\frac{c}{1-\epsilon}\cdot t\cdot \varphi^*(t)\
\end{aligned}
\end{equation*}
\begin{equation*}
\begin{aligned}
 \Rightarrow \frac{1}{t} \int_{0}^{t}\varphi^*(y)dy\leq \frac{c}{1-\epsilon}\varphi^*(t)
\phantom{\hspace{7cm}},
\end{aligned}
\end{equation*}
for every $\epsilon>0$. Letting $\epsilon\rightarrow 0^+$, we conclude \eqref{eq:1p2} for any $t\in(0,1]$.

\section{$L^p$ integrability of $A_1$ weights on $\mathbb{R}.$} \label{sec:3}
We will now prove the following:
\begin{lemma} \label{lem:3p1}
Let $g:(0,1)\rightarrow \mathbb{R}^+$ be a non-increasing, left continuous function which satisfies the following inequality:
\begin{equation}
\frac{1}{t} \int_{0}^{t}g(y)dy\leq c\cdot g(t), \forall t\in (0,1]
\end{equation}
for a fixed $c>1$. Then for any $p\in [1,\frac{c}{c-1})$ the following is the:
\begin{equation}
\int_{0}^{1}g^p(y)dy\leq \frac{1}{c^{p-1}(c+p-pc)}\Big( \int_{0}^{1}g(y)dy\Big)^p.
\end{equation}
Moreover, inequality (3.2) is best possible.
\end{lemma}

$\underline{\textbf{Proof}}:$
Fix a $p$ such that $1\leq p< \frac{c}{c-1}$ and let $F=\int_{0}^{1}g^p(y)dy$ and $f=\int_{0}^{1}g(y)dy$. Then by H\"{o}lder's inequality $f^p\leq F$. We need to prove that
\begin{equation}
 F\leq \frac{1}{c^{p-1}(c+p-pc)}\cdot f^p.
\end{equation}
We define the following function:
\begin{equation*}
H_p:\Big[ 1,\frac{p}{p-1}\Big] \rightarrow [0,1]
\end{equation*}
by $H_p(z)=pz^{p-1}-(p-1)z^p$. Then we easily see that $H_p$ is one to one and onto. We denote it's inverse function by $\omega_p$ defined on $[0,1]$, which is decreasing as $H_p$ also is. We shall prove that (3.3) holds, equivalently, $H_p(c)\leq \frac{f^p}{F}$ $\Leftrightarrow c\geq \omega_p \Big(\frac{f^p}{F}\Big)=: \tau$.
\newline
Suppose on the contrary that $c<\tau$. We are going to reach to a contradiction.

Define the following function on $(0,1]$, by $g_1(t)=\frac{f}{\tau}t^{-1+\frac{1}{\tau}}$. This is obviously non-increasing and continuous $(0,1]$. Additionally, it satisfies for any $t\in (0,1]$ the following equality.
 \begin{equation}
 \frac{1}{t}\int_{0}^{t}g_1(y)dy=\tau \cdot g_1(t)
\end{equation}
Indeed: $\frac{1}{t} \int_{0}^{t}g_1(y)dy=\frac{1}{t} \frac{f}{\tau}\int_{0}^{t}y^{-1+\frac{1}{\tau}}dy=\frac{f}{t}\big[y^{\frac{1}{\tau}}\big]_{y=0}^t=\frac{f}{t}\cdot t^{\frac{1}{\tau}}=\tau \cdot \Big(\frac{f}{\tau}t^{-1+\frac{1}{\tau}}\Big)=\tau g_1(t)$. Moreover, it satisfies $\int_{0}^{1}g(y)dy=f$ and $\int_{0}^{1}g^p(y)dy=F$. The first equation is obvious, in view of (3.4). As for the second it is equivalent to
$\frac{f^p}{\tau^p}\int_{0}^{1}y^{-p+\frac{p}{\tau}}dy=F$  $\Leftrightarrow \frac{f^p}{\tau^p(1+\frac{p}{\tau}-p)}=F$ $\Leftrightarrow p\tau^{p-1}-(p-1)\tau^p=\frac{f^p}{F}$ $\Leftrightarrow H_p(\tau)=\frac{f^p}{F}$ $\Leftrightarrow \tau=\omega_p(\frac{f^p}{F})$, which is true by the definition of $\tau$.
%\newline

We are now aiming to prove that the following inequality is satisfied:
\begin{equation}
\int_{0}^{t}g(y)dy\leq \int_{0}^{t} g_1(y)dy, \text{for any t}\in (0,1].
\end{equation}
For this reason we define the following subset of $(0,1)$:
\newline
$G=\Big\{t\in (0,1): \int_{0}^{t}g(y)dy> \int_{0}^{t}g_1(y)dy\Big\}$, and we suppose that G is non empty. By the continuity of the involving integral functions on $t$  we have as a consequence that G is an open subset of $(0,1)$. Since $G\neq \emptyset$ $\Rightarrow G=\bigcup\limits_\nu I_\nu$, where $(I_\nu)_\nu$ is a (possibly finite) sequence of pairwise disjoint open integrals on $(0,1)$. Let us choose one of them, $I_\nu=(\alpha_\nu,b_\nu)$. Since $\alpha_\nu \notin G$
\begin{equation}
\Rightarrow \int_{0}^{\alpha_\nu} g(y)dy \leq \int_{0}^{\alpha_\nu} g_1(y)dy.
\end{equation}
Let now $(x_n)_n\subseteq I_\nu$ be a sequence such that $x_n\rightarrow \alpha_\nu$, as $n\rightarrow \infty$. Since $x_n\in G, \forall n=1,2,...$ we must have that $\int_{0}^{x_n}g(y)dy> \int_{0}^{x_n}g_1(y)dy$, so letting $n\rightarrow \infty$ we conclude that
\begin{equation}
\int_{0}^{\alpha_\nu}g(y)dy\geq \int_{0}^{\alpha_\nu}g_1(y)dy.
\end{equation}
By (3.6) and (3.7) we see that $\int_{0}^{\alpha_\nu}g(y)dy=\int_{0}^{\alpha_\nu}g_1(y)dy$. In the same way we prove that $\int_{0}^{b_\nu}g(y)dy=\int_{0}^{b_\nu}g_1(y)dy$. As a consequence, we must have that
\begin{equation}
\int_{\alpha_\nu}^{b_\nu}g(y)dy=\int_{\alpha_\nu}^{b_\nu}g_1(y)dy.
\end{equation}
Let now $t\in I_\nu=(\alpha_\nu,b_\nu)$. Since $t\in G$ and because  of (3.1) and (3.4) and the assumption on $\tau$, we must have the following:
$cg(t)\geq \frac{1}{t}\int_{0}^{t}g(y)dy>\frac{1}{t}\int_{0}^{t}g_1(y)dy=\tau \cdot g_1(t)>c g_1(t)$ thus $g(t)>g_1(t)$, for every $t\in I_\nu$. This is impossible in view of (3.8).
We note the following (which can be seen in \cite{5}, p.88).

\begin{lemma} \label{lem:3p1}
Let $\varphi_1, \varphi_2: (0,1]\rightarrow \mathbb{R}^+$ be integrable functions. Then the following are equivalent
\begin{enumerate}[i)]{\em
\item $\int_{0}^{t}\varphi_1^*(y)dy \leq \int_{0}^{t}\varphi_2^*(y)dy$, for every $t\in (0,1]$.
\item $\int_{0}^{1}G(\varphi_1(x))dx  \leq \int_{0}^{1}G(\varphi_2(x))dx$
\vspace{-5pt}}
\end{enumerate}
for any G convex, non-negative, increasing and left continuous function on $[0,+\infty)$.
\end{lemma}
%\newline
We consider now two cases:
\begin{enumerate}[A)]
\item We have equality in $(3.5)$ for every $t\in (0,1]$. That is $\int_{0}^{t}g(y)dy =\int_{0}^{t}g_1(y)dy$, $ \forall t\in (0,1]$. This gives immediately as a consequence that $g(t)=g_1(t)$ almost everywhere on $(0,1]$, and since $g_1$ is continuous on $(0,1]$, we must have that $g(t)=g_1(t), \forall t\in (0,1]$ $\Rightarrow g(t)=\frac{f}{\tau}t^{-1+\frac{1}{\tau}}, \forall t\in (0,1]$ $\Rightarrow \frac{1}{t}\int_{0}^{t}g(y)dy=\tau g(t), \forall t\in (0,1]$.
Then in view of $(3.1)$ we conclude that $c\geq \tau$ which is a contradiction since we have supposed the opposite inequality.
\item There exists a $t_0 \in (0,1)$ such that:
\begin{equation*}
\int_{0}^{t_0}g(y)dy<\int_{0}^{t_0}g_1(y)dy.
\end{equation*}
Then, by continuity reasons, we have as a consequence that there exists a $\delta>0$ such that
\begin{equation}
\int_{0}^{t_0}g(y)dy<\int_{0}^{t_0}g_1(y)dy, \,\text{for any t} \in(t_0-\delta, t_0+\delta)=I_\delta.
\end{equation}
\vspace{-5pt}
\end{enumerate}
We define now the quantities $d_1,d_2$ by the following equations:
\begin{equation}
\frac{1}{\delta}\int_{t_0-\delta}^{t_0}g_1(y)dy=d_1 \quad \text{and} \quad \frac{1}{\delta}\int_{t_0}^{t_0+\delta}g_1(y)dy=d_2.
\end{equation}
Then by H\"{o}lder's inequality on the interval $(t_0-\delta,t_0)$ for $g_1$, we conclude that
\begin{equation}
\frac{1}{\delta}\int_{t_0-\delta}^{t_0}g_1^p(y)dy>d_1^p,
\end{equation}
which is a strict inequality since $g_1$ is strictly decreasing (therefore not constant) on the interval $(t_0-\delta,t_0)$. In the same way we have
\begin{equation}
\frac{1}{\delta}\int_{t_0}^{t_0+\delta}g_1^p(y)dy>d_2^p.
\end{equation}
Then since $g_1$ is decreasing we have that $d_2<d_1$.
We define now the following nonincreasing (as can be easily seen) function on $(0,1]$
\begin{equation}
 g_2(t) = \left\{ \begin{array}{ll}
         g_1(t), & \mbox{$t\in (0,1]\setminus (t_0-\delta,t_0+\delta)$}\\
        d_1 & \mbox{$t\in[t_0-\delta,t_0)$}\\
        d_2 & \mbox{$t\in[t_0,t_0+\delta]$}.\end{array} \right.
\end{equation}
By $(3.9)$ and since $g_1$ is decreasing we easily see that we can choose $\delta>0$ small enough, so that
\begin{equation}
\int_{0}^{t}g(y)dy\leq \int_{0}^{t}g_2(y)dy, \,\text{for any}\, t\in(0,1].
\end{equation}
Additionally, because of $(3.11)$ and $(3.12)$ we must have that
\begin{equation*}
\int_{0}^{1}g_2^p(y)dy< \int_{0}^{1}g_1^p(y)dy=F.
\end{equation*}
Since $(3.14)$ holds for any $t\in(0,1]$ and because of Lemma 3.2 we conclude that $\int_{0}^{1}g^p(y)dy\leq \int_{0}^{1}g_2^p(y)dy<F$, by considering the function $G(t)=t^p$. This is obviously a contradiction according to the way that F is defined. In this way we derive the proof of our Lemma.
\newline
We now proceed to the:
\newline
$\underline{\textbf{Proof of Theorem 2}}:$ Without loss of generality we suppose that $\mc J=(0,1)$. Let $p\in[1,\frac{c}{c-1})$ and $\mc I\subseteq (0,1)$ and let also $\varphi_{\mc I}=\varphi/_{\mc I}$ the restriction of $\varphi$ to $\mc I$. Consider now the function $g: (0,|\mc I|] \rightarrow \mathbb{R}^+$, defined by $g=(\varphi_{\mc I})^*$. Then since $\varphi_{\mc I}\in A_1(\mc I)$ with $A_1$ constant not more than $c$, we must have by using Theorem 1 that $\frac{1}{t}\int_{0}^{t}g(y)dy\leq cg(t)$, for any $t\in (0,|\mc I|]$.
Thus by Lemma (3.1) we have as a consequence that
\begin{equation*}
\frac{1}{|\mc I|}\int_{0}^{|\mc I|}g^p(y)dy\leq \frac{1}{c^{p-1}(c+p-pc)} \Big(\frac{1}{|\mc I|} \int_{0}^{|\mc I|}g(y)dy \Big)^p,
\end{equation*}
which is:
\begin{equation*}
\frac{1}{|\mc I|}\int_{\mc I}\varphi^p(x)dx \leq \frac{1}{c^{p-1}(c+p-pc)} \Big(\frac{1}{|\mc I|} \int_{\mc I}\varphi(x)dx \Big)^p.
\end{equation*}
The relation (1.3) is proved.

Nikolidakis Eleftherios, Visiting Professor, University of Ioannina, Department of Mathematics, GR 45110, Panepistimioupolis, Greece.

\end{document}